# Insight into Primal Augmented Lagrangian Multilplier Method

Premjith B[1], Sachin Kumar S[2], Akhil Manikkoth, Bijeesh T V, K P Soman
Centre for Excellence in Computational Engineering and Networking
Amrita Vishwa Vidyapeetham
Coimbatore, India
[1]Email: prem.jb@gmail.com
[2]Email: sachinnme@gmail.com

*Abstract*—We provide a simplified form of Primal Augmented Lagrange Multiplier algorithm. We intend to fill the gap in the steps involved in the mathematical derivations of the algorithm so that an insight into the algorithm is made. The experiment is focused to show the reconstruction done using this algorithm.

## I. INTRODUCTION

Compressive Sensing (CS) is one of the hot topics in the area of signal processing [1,2,3]. The conventional way to sample the signals follows the Shannon's theorem, i.e., the sampling rate must be at least twice the maximum frequency present in the signal (Nyquist rate) [4]. For practical signals using CS, the sampling or sensing goes against the Nyquist rate. Consider the sensing matrix, $A$ ($A \in \mathbb{R}^{nxN}$), as the concatenation of two orthogonal matrices, $\psi$ and $\phi$. A signal, $b \in \mathbb{R}^n$ (non-zero vector) can be represented as the linear combination of columns of $\Psi$ or as a linear combination of columns of $\phi$. That is in and $b = \phi\beta$, $\alpha$ and $\beta$ are defined clearly. Taking $\psi$ as identity matrix and $\phi$ as Fourier transform matrix, we can infer that α is the time domain representation of b and $\beta$ is the frequency-domain representation. Here we can see that a signal cannot be sparse in both domains at a time. So this phenomena can be extended as, for any arbitrary pair of bases $\psi$ and $\phi$, either $\alpha$ can be sparse or $\beta$ can be sparse. Now what if the bases are same, the vector b can be constructed using the one of the columns in $\psi$ and get the smallest possible cardinality (smallest possible number) in $\alpha$ and $\beta$. The proximity between the two bases can be defined through mutual-coherence which can be defined as the maximal inner product between columns from two bases. [5,6]. The CS theory says that the signal can be recovered from very few samples or measurements. This is made true based on two principles: *sparsity* and *incoherence.*

*Sparsity*: a signal is said to be sparse if it is represented using much less number of sample coefficients without loss of information. The advantage is that sparsity gives fast calculation. CS exploits the fact that the natural signals are sparse in nature when expressed using proper basis [5, 6].

*Incoherence:* we can understand mutual-coherence between the measurement basis and sparsity basis as

$$\mu = \sqrt{N} \max_{i,j} |<\phi_i, \psi_j>|.$$

For matrices which are having orthonormal basis, $\mu = 1$. If the incoherence is higher, the number of measurements required will be smaller [7].

### A. Underdetermined Linear System:

Consider a matrix $A \in \mathbb{R}^{nxN}$ with $n < N$, and define the underdetermined linear system of equations $Ax = b$, $x \in \mathbb{R}^N$, $b \in \mathbb{R}^n$. This means that there are many ways of representing a given signal $b$ (this system has more unknowns than equations, and thus it has either no solution, if b is not in the span of the columns of the matrix $A$, or infinitely many solutions). In order to avoid the anomaly of having no solution, we shall hereafter assume that $A$ is a full-rank matrix, implying that its columns span the entire space $\mathbb{R}^n$. Consider an example of image reconstruction. Let $b$ be an image with lower quality, and we need to reconstruct the original image, which is represented as $x$ using a matrix $A$. Recovering $x$ given $A$ and $b$ constitutes the linear inversion problem. Compressive sensing is the method to find the sparsest solution to some underdetermined system of linear equations. In compressive sensing point of view $b$ is called the measurement vector, $A$ is called the sensing matrix or measurement matrix and $x$ is the unknown signal. A conventional solution is using the linear least square method which increases the computational time.[8].

The above mentioned problem can be solved through optimization method using $l_0$-minimization.

$$\min_x \|x\|_0 \text{ Subject to } Ax = b$$

The solution to this problem is to find an $x$ that will have few non-zero entries. This is known as $l_0$ norm.

The properties that a norm should satisfy are

(1) Zero vectors: $\|v\| = 0$ if and only if $\forall v = 0$

(2) Absolute homogeneity:, $\forall t \neq 0$
$\|tu\| = \|t\|\|u\|, t \in \mathbb{R}$

(3) Triangle inequality: $\|u + v\| \leq \|u\| + \|v\|$

This problem is NP-Hard. Candes and Donoho have proved if the signal is sufficiently sparse, the solution can be obtained using $l_1$-minimization. The objective function becomes,

$$\min_x \|x\|_1 \text{ Subject to } Ax = b$$

Which is a convex optimization problem. This minimization tries to find out a vector $x$ whose absolute sum is small among all the possible solutions. This is known as $l_1$ norm. This problem can be recasted as the linear problem and can be tried to solve using interior-point method. But this becomes computationally complex. Most of the real world applications need to robust to noise, that is, the observation vector, $b$, can be corrupted by noise. To handle the noise a slight change in the constraint is made. An additional threshold parameter is added which is predetermined based on the noise level.

$$\min_x \|x\|_1 \text{ Subject to } \|b - Ax\|_2 \leq T$$

$T > 0$ is the threshold. For the reconstruction of such signals basis pursuit denoising (BPDN) algorithm is used. [5, 9, 10].

In the light of finding efficient algorithms to solve these problems, several algorithms have been proposed like, Orthogonal Matching Pursuit, Primal-Dual Interior-Point Method, Gradient Projection, Homotopy, Polytope Faces Pursuit, Iterative Thresholding, Proximal Gradient, Primal Augmented Lagrange Multiplier, Dual Augmented Lagrange Multiplier [11]. These algorithms can work better when the signal representation is more sparse. Clearly there are infinitely many solutions out which few will give the good reconstruction. To find out a single solution, a regularization in used. Finally what is trying to achieve is a solution which has the minimum norm.

## II. PRIMAL AUGMENTED LAGRANGE MULTIPLIER METHOD

Our aim is to solve the system $Ax = b$. But in practical case there could be some error. So we write $Ax \approx b$. So when it is written in equality form, it ill be $Ax + r = b$, where $r$ is the error of residual. So our main objective is to minimize $x$ and $r$, such that $x$ is as sparse as possible. So we take $l_1$ norm of $x$.

The objective function is,

$$\min_{x,r} \|x\|_1 + \frac{1}{2\mu}\|r\|^2 \quad (1)$$

Subject to, $Ax + r = b$

Where $X \in C^n$, $r \in C^m$

Take the Lagrangian of the form,

$$L(x, r, y) = \min_{x,r} \|x\|_1 + \frac{1}{2\mu}\|r\|^2 - \text{Re}(y*(Ax + r - b)) +$$

$$+ \frac{\beta}{2}\|Ax + r - b\|^2 \quad (2)$$

Where, $y \in C^m$ is a multiplier and $\beta > 0$ is a penalty parameter. Since $y \in C^m$, instead of $y^T$ we use $y*$ (conjugate transpose).

Using iterative methods, for a given $(x^k, y^k)$ we obtain $(r^{k+1}, x^{k+1}, y^{k+1})$. We use alternating minimization method to achieve this. *i.e.,* for fixed values of two variables we minimize the other. First, fix $x = x^k$ and $y = y^k$ and find $r^{k+1}$.

Minimize $L(x, r, y)$ with respect to $r$

Find the derivative of $L(x, r, \mu)$ with respect to $r$

$$\frac{\partial L(x,r,y)}{\partial r} = \frac{\partial \|x\|_1}{\partial r} + \frac{1}{2\mu} \frac{\partial \|r\|^2}{\partial r}$$

$$- \frac{\partial \operatorname{Re}(y * (Ax + r - b))}{\partial r} + \frac{\beta}{2} \frac{\partial \|Ax + r - b\|^2}{\partial r} \quad (3)$$

Omit all the terms independent of $r$.

$\|r\|^2$ Can be considered as $r^T r$

So, eq. (3) will be,

$$\frac{\partial L(x,r,\mu)}{\partial r} = \frac{1}{2\mu} \frac{\partial r^T r}{\partial r} - \frac{\partial \operatorname{Re}(y * (Ax + r - b))}{\partial r} +$$

$$\frac{\beta}{2} \frac{\partial \|Ax + r - b\|^2}{\partial r} = 0$$

$$\Rightarrow \frac{2r}{2\mu} - \frac{\partial (yr^T)}{\partial r} +$$

$$\frac{\beta}{2} \frac{\partial ((Ax + r) - b)^T ((Ax + r) - b)}{\partial r} = 0$$

$$\Rightarrow \frac{r}{\mu} - y + \frac{\beta}{2} \frac{\partial ((Ax + r)^T (Ax + r) - 2b(Ax + r)^T)}{\partial r} = 0$$

$$\Rightarrow \frac{r}{\mu} - y + \frac{\beta}{2} \frac{\partial ((Ax)^T (Ax) + 2r^T (Ax) + r^T r - 2b(Ax + r)^T)}{\partial r} = 0$$

$$\Rightarrow \frac{r}{\mu} - y + \frac{\beta}{2} (2Ax + 2r - 2b) = 0$$

$$\Rightarrow \frac{r}{\mu} - y + \beta(Ax + r - b) = 0$$

$$\Rightarrow \frac{r}{\mu} + \beta r - y + \beta(Ax - b) = 0$$

$$\Rightarrow r\left(\frac{1}{\mu} + \beta\right) = y - \beta(Ax - b)$$

$$\Rightarrow r\left(\frac{1 + \mu\beta}{\mu}\right) = y - \beta(Ax - b)$$

$$\Rightarrow r\left(\frac{1 + \mu\beta}{\mu\beta}\right) = \frac{y}{\beta} - (Ax - b)$$

$$\Rightarrow r = \left(\frac{\mu\beta}{1 + \mu\beta}\right)\left(\frac{y}{\beta} - (Ax - b)\right)$$

$$\Rightarrow r^{k+1} = \left(\frac{\mu\beta}{1 + \mu\beta}\right)\left(\frac{y^k}{\beta} - (A^x - b)\right)$$

(4)

Now, fix $r = r^{k+1}$ and $y = y^k$. Do minimization of the Lagrangian with respect to $x$.

$$\min_x \|x\|_1 + \frac{1}{2\mu} \|r\|^2 - \operatorname{Re}(y * (Ax + r - b))$$

$$+ \frac{\beta}{2} \|Ax + r - b\|^2$$

(5)

Omit all the terms independent of $x$.

$$\Rightarrow \min_x \|x\|_1 - \operatorname{Re}(y * (Ax)) +$$

$$\frac{\beta}{2} (Ax + r - b)^T (Ax + r - b)$$

$$= \min_x \|x\|_1 - \operatorname{Re}(y * (Ax)) +$$

$$\frac{\beta}{2} ((Ax + r)^T (Ax + r) - 2b^T (Ax + r))$$

$$= \min_{x} \|x\|_1 - \text{Re}(y * (Ax)) + \frac{\beta}{2}((Ax)^T (Ax) + 2r^T Ax - 2b^T (Ax + r))$$

$$= \min_{x} \|x\|_1 - \text{Re}(y * (Ax)) + \frac{\beta}{2}(x^T A^T Ax + 2r^T Ax - 2b^T Ax)$$

Here we assume $A^T A = I$

$$\Rightarrow \min_{x} \|x\|_1 + \frac{\beta}{2}\left(x^T x + 2r^T Ax - 2b^T Ax - \frac{2y^T Ax}{\beta}\right) \quad (6)$$

This can be written as,

$$\min_{x} \|x\|_1 + \frac{\beta}{2}\left\|Ax + r^{k+1} - b - \frac{y^k}{\beta}\right\|^2 \quad (7)$$

Constant terms in the expansion of $\left\|Ax + r^{k+1} - b - \frac{y^k}{\beta}\right\|^2$ make no sense.

Since (7) is in quadratic form, we can approximate it by first two terms of Taylor series.

$$\min_{x} \|x\|_1 + \beta\left(\text{Re}((g^k) * (x - x^k))\right) + \frac{1}{2\tau}\|x - x^k\|^2 \quad (8)$$

Where $\tau > 0$ is a proximal parameter and,

$$g^k \overset{\Delta}{=} A * \left(Ax^k + x^{k+1} - b - \frac{y^k}{\beta}\right) \quad (9)$$

$g^k$ is the gradient of eq. (7) at $x = x^k$ excluding the multiplication by $\beta$. The solution of eq. (8) is given by,

$$x^{k+1} = \text{Shrink}\left(x^k - \tau g^k, \frac{\tau}{\beta}\right) \quad (10)$$

Finally update the multiplier $y$ by fixing, $r = r^{k+1}$ and $x = x^{k+1}$

Since the Lagrangian is quadratic in $y$ and for any fixed $x^{k+1}$ and $r^{k+1}$,

$$y^{k+1} = y^k - \gamma\beta\left(Ax^{k+1} + r^{k+1} - b\right) \quad (11)$$

Where $\gamma > 0$ is a constant.

So the PALM algorithm in 3 steps is,

1. $r^{k+1} = \left(\frac{\mu\beta}{1+\mu\beta}\right)\left(\frac{y^k}{\beta} - (A^x - b^T)\right)$

2. $x^{k+1} = \text{Shrink}\left(x^k - \tau g^k, \frac{\tau}{\beta}\right)$

3. $y^{k+1} = y^k - \gamma\beta\left(Ax^{k+1} + r^{k+1} - b\right)$

III. RESULTS

Primal Augmented Lagrangian Multiplier algorithm is tested on three different images of size 256x256 pixels. Then Gaussian, Salt & pepper and Speckle noises with varying percentage of intensities are added onto these images and reconstruct the image using PALM algorithm. PSNR, RMSE and execution time are calculated in each experiment and reported in the tables Table1 and Table2.

In the case of first two images, PSNR is less and RMSE is more for original image. But in the case of third image, when we add 20% Gaussian noise, it shows good recovery.

## IV. CONCLUSION

We have presented a simple mathematical derivation for Primal Augmented Lagrangian Multiplier method which is considered as one of the best compressive sensing algorithms. PALM can be used for reconstructing various signals with varying noise intensities. Experiments show that this algorithm gives very good reconstruction result.

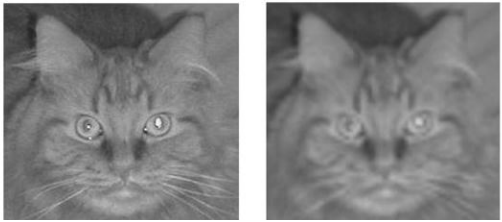

Figure 1: (a) Original Image (b) Reconstructed Image

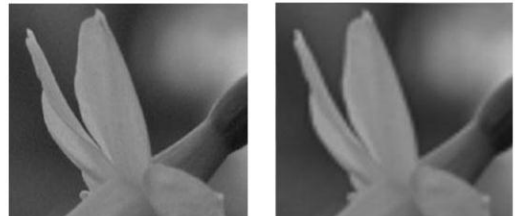

Figure 1: (a) Original Image (b) Reconstructed Image

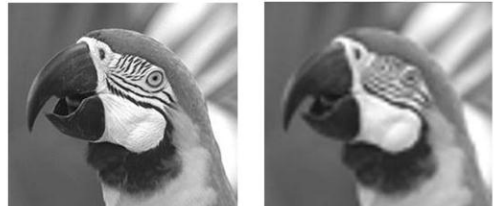

Figure 3: (a) Original Image (b) Reconstructed Image

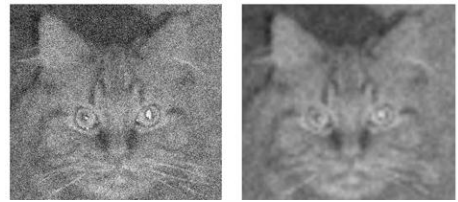

Figure 3: Added with 2% Gaussian noise

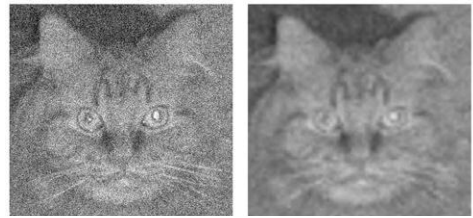

Figure 4: Added with 5% Gaussian noise

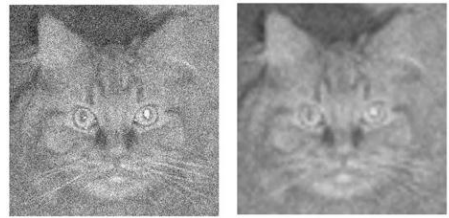

Figure 5: Added with 10% Gaussian noise

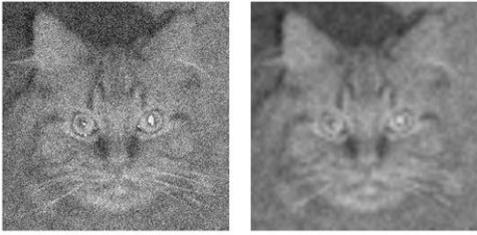

Figure 6: Added with 20% Gaussian noise

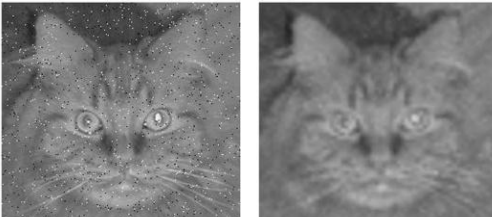

Figure 7: Added with 2% Salt & pepper noise

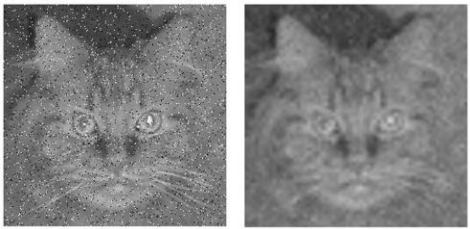

Figure 8: Added with 5% Salt & pepper noise

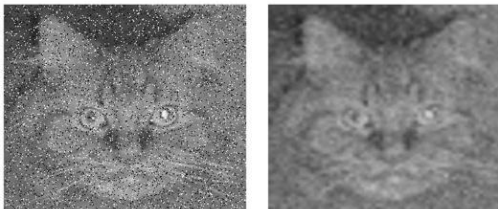

Figure 9: Added with 10% Salt & pepper noise

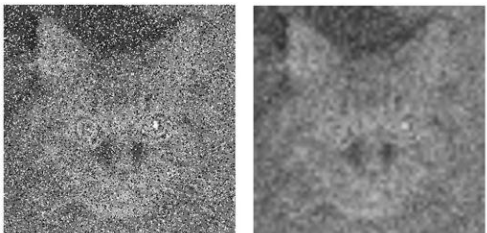

Figure 10: Added with 20% Salt & pepper noise

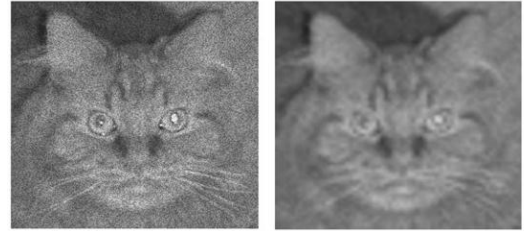

Figure 11: Added with 2% Speckle noise

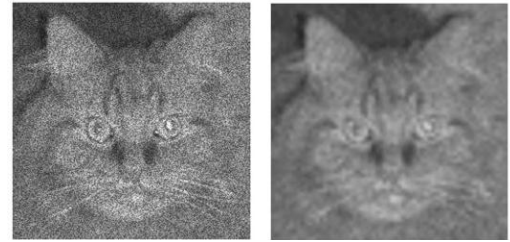

Figure 12: Added with 5% Speckle noise

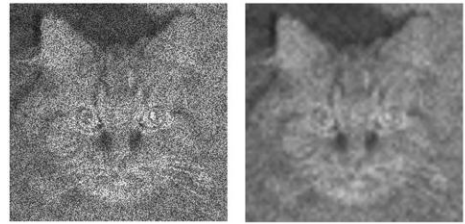

Figure 13: Added with 10% Speckle noise

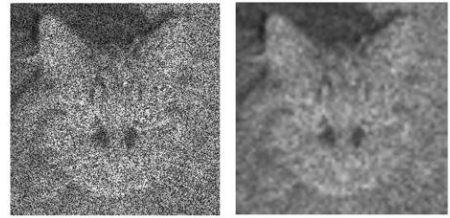

Figure 14: Added with 20% Speckle noise

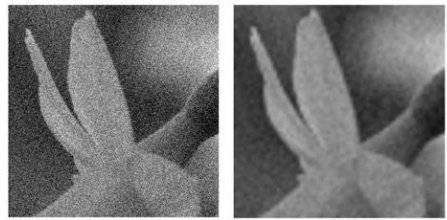

Figure 15: Added with 2% Gaussian noise

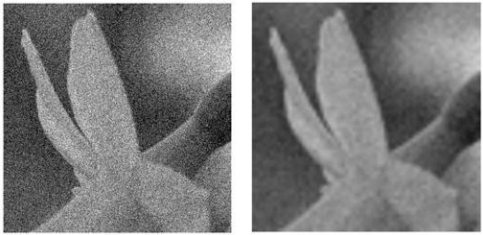

Figure 16: Added with 5% Gaussian noise

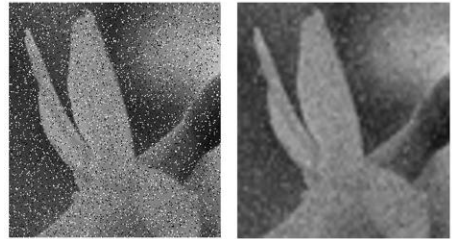

Figure 21: Added with 10% Salt & pepper noise

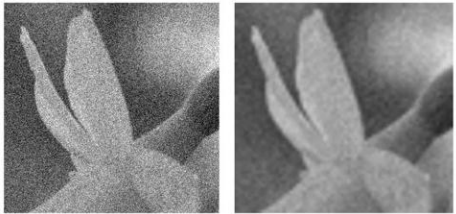

Figure 17: Added with 10% Gaussian noise

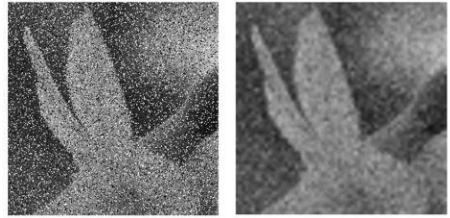

Figure 22: Added with 20% Salt & pepper noise

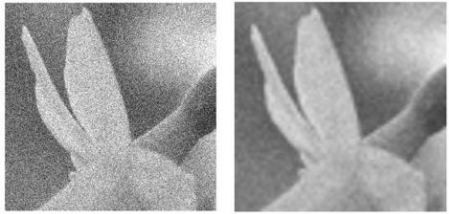

Figure 18: Added with 20% Gaussian noise

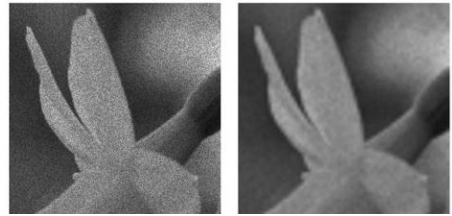

Figure 23: Added with 2% Speckle noise

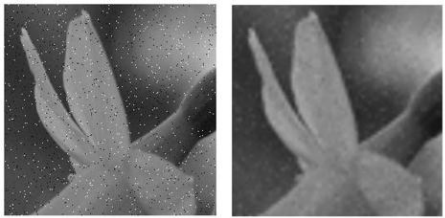

Figure 19: Added with 2% Salt & pepper noise

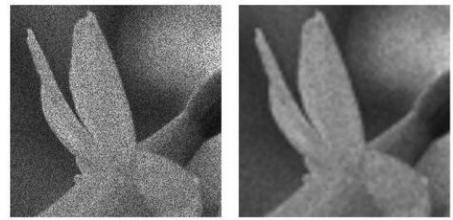

Figure 24: Added with 5% Speckle noise

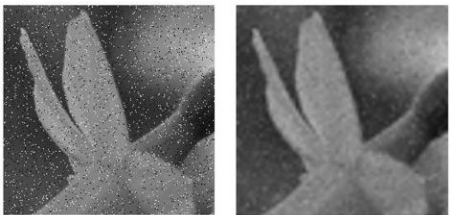

Figure 20: Added with 5% Salt & pepper noise

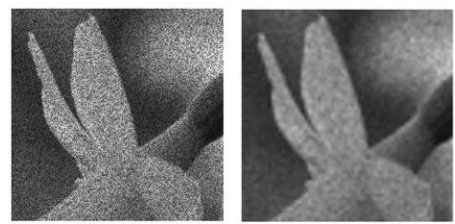

Figure 25: Added with 10% Speckle noise

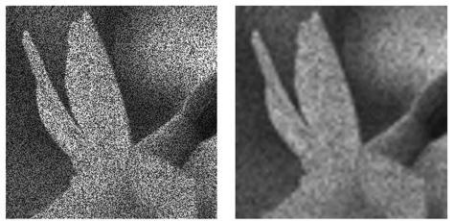

Figure 26: Added with 20% Speckle noise

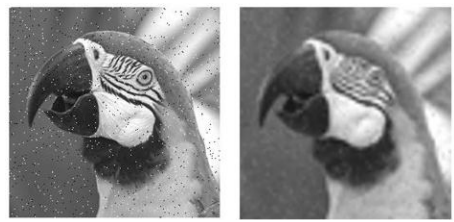

Figure 31: Added with 2% Salt & pepper noise

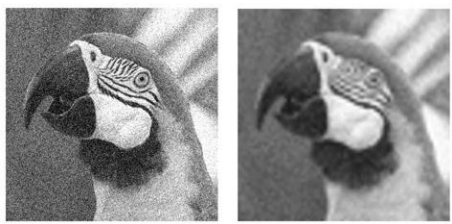

Figure 27: Added with 2% Gaussian noise

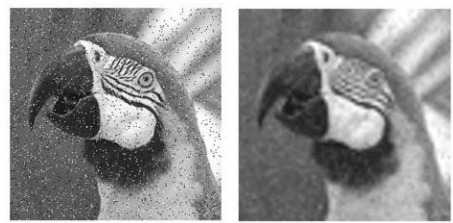

Figure 32: Added with 5% Salt & pepper noise

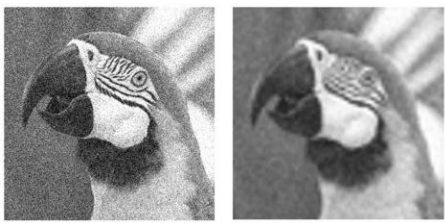

Figure 28: Added with 5% Gaussian noise

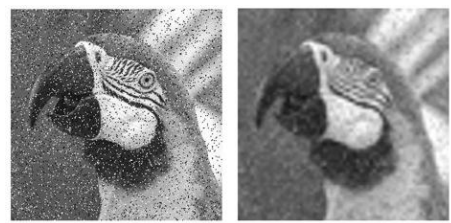

Figure 33: Added with 10% Salt & pepper noise

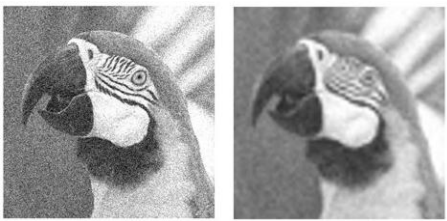

Figure 29: Added with 10% Gaussian noise

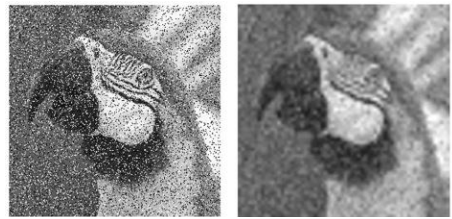

Figure 34: Added with 20% Salt & pepper noise

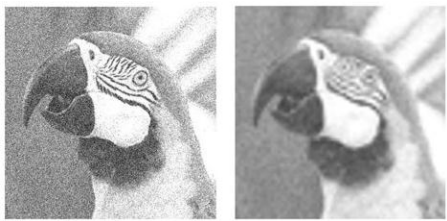

Figure 30: Added with 20% Gaussian noise

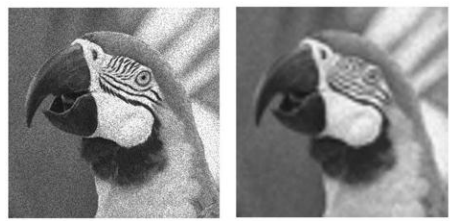

Figure 35: Added with 2% Speckle noise

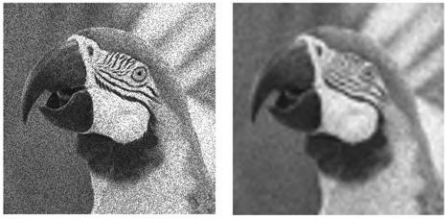

Figure 36: Added with 5% Speckle noise

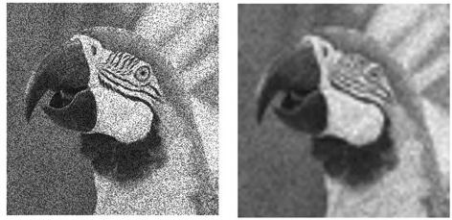

Figure 37: Added with 10% Speckle noise

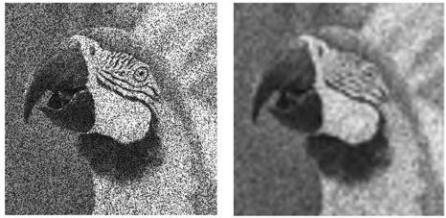

Figure 38: Added with 20% Speckle noise

TABLE 1: PSNR, RMSE AND EXECUTION TIME OF ORIGINAL IMAGES

| Image name | PSNR (dB) | RMSE | Time (Sec) |
|---|---|---|---|
| Image_1 | 17.8557 | 32.6403 | 20.940833 |
| Image_2 | 14.0948 | 50.3272 | 19.680670 |
| Image_3 | 10.5039 | 76.0933 | 21.732722 |

TABLE 2: PSNR, RMSE AND EXECUTION TIME OF ORIGINAL IMAGES WHEN ADDED WITH NOISES OF VARYING INTENSITIES

| Image name | Noise | % of error | PSNR (dB) | RMSE | Time (Sec) |
|---|---|---|---|---|---|
| Image_1 | Gaussian | 2 | 15.721 | 41.734 | 23.743473 |
| | | 5 | 15.77 | 41.4995 | **19.613224** |
| | | 10 | 15.7262 | 41.7089 | 25.044838 |
| | | 20 | **15.7814** | **41.445** | 26.502142 |
| | Salt & pepper | 2 | **16.8149** | **36.7956** | 28.512797 |
| | | 5 | 15.5122 | 42.7494 | **22.795465** |
| | | 10 | 14.009 | 50.8263 | 27.19995 |
| | | 20 | 11.9891 | 64.1333 | 27.10357 |
| | Speckle | 2 | **16.7394** | **37.1169** | 23.400962 |
| | | 5 | 15.5125 | 42.7481 | 22.950852 |
| | | 10 | 13.987 | 50.9556 | **22.495303** |
| | | 20 | 12.0597 | 63.6144 | 23.600862 |
| Image_2 | Gaussian | 2 | 13.0878 | 56.5134 | 28.182993 |
| | | 5 | **13.1057** | **56.3966** | 24.715946 |
| | | 10 | 13.0987 | 56.4423 | **24.711972** |
| | | 20 | 13.0964 | 56.4571 | 26.883013 |
| | Salt & pepper | 2 | **13.6543** | **52.9451** | 26.492682 |
| | | 5 | 13.0383 | 56.8364 | **21.902348** |
| | | 10 | 12.2250 | 62.4152 | 23.998090 |
| | | 20 | 10.8429 | 73.1803 | 26.126452 |
| | Speckle | 2 | **13.7082** | **52.6175** | 28.975433 |
| | | 5 | 13.1603 | 56.0435 | **25.038859** |
| | | 10 | 12.3823 | 61.2953 | 28.131446 |
| | | 20 | 11.3225 | 69.2489 | 25.585058 |
| Image_3 | Gaussian | 2 | 10.2427 | 78.4158 | 23.582806 |
| | | 5 | 10.3201 | 77.7214 | **21.334992** |
| | | 10 | 10.5091 | 76.0479 | 26.421035 |
| | | 20 | **11.2182** | **70.0861** | 21.982944 |
| | Salt & pepper | 2 | **10.3787** | **77.1977** | 21.228414 |
| | | 5 | 10.1755 | 79.0254 | 21.854468 |
| | | 10 | 9.8367 | 82.1689 | 25.088871 |
| | | 20 | 9.3083 | 87.3222 | 25.632787 |
| | Speckle | 2 | **10.4489** | **76.5762** | 21.490942 |
| | | 5 | 10.3024 | 77.8793 | 22.401311 |
| | | 10 | 10.0795 | 79.9031 | 24.353466 |
| | | 20 | 9.6208 | 84.2366 | 24.639351 |